\documentstyle[12pt,twoside]{amsart}

\newcommand{\p}[0]{{\mathbb P}}

\newcommand{\map}[0]{\dasharrow}

\newcommand{\chr}[0]{\operatorname{char}}

\newsymbol\subsetneq 2328
\newsymbol\onto 1310

\newsymbol\twoheadrightarrow 1310

\newtheorem{thm}{Theorem}

\newtheorem{lem}[thm]{Lemma}

\newtheorem{prop}[thm]{Proposition}

\theoremstyle{definition}

\newtheorem{rem}[thm]{Remark}          

\newtheorem{defn-lem}[thm]{Definition--Lemma}

\theoremstyle{remark}

\begin{document}
\bibliographystyle{amsplain}

\title{Fixed points of group actions and rational maps}
\author{J{\'a}nos Koll{\'a}r and Endre Szab{\'o}}

\maketitle

The aim of this note is to give simple proofs of
the results in section 5 of
\cite{ry} about the behaviour of fixed points of finite group
actions under rational maps. Our proofs   work in any
characteristic.

\begin{lem}\label{class.defn}
 Let $K$ be an algebraically closed field and  $H$
a (not necessarily connected) linear
algebraic group over $K$.  The following are equivalent.
\begin{enumerate}
\item Every representation  $H\to GL(n,K)$ has an
$H$-eigenvector.
\item There is a  (not necessarily connected) unipotent, normal 
subgroup $U<H$ such that $H/U$ is abelian.
\end{enumerate}
\end{lem}

Proof. Let $H\to GL(n,K)$ be a faithful representation.
If (\ref{class.defn}.1) holds then $H$ is conjugate to
an upper triangular subgroup, this implies (\ref{class.defn}.2).

Conversely, any representation of a unipotent group
has eigenvectors (with eigenvalue 1, cf.\ \cite[I.4.8]{borel})
and the subspace of all eigenvectors is an
$H/U$-representation.\qed

\begin{prop}[Going down]\label{Going-down}
  Let $K$ be an algebraically closed
field,
$H$ a linear  algebraic group over $K$ and $f:X\map Y$ an
$H$-equivariant map of
$K$-schemes. Assume that
\begin{enumerate}
\item $H$ satisfies the equivalent conditions
of (\ref{class.defn}),
\item $H$ has a smooth fixed point on $X$, and
\item $Y$ is proper.
\end{enumerate}
\noindent Then $H$ has a fixed point on $Y$.
\end{prop}

Proof. The proof is by induction on $\dim X$. The  case  $\dim
X=0$ is clear.

Let $x\in X$ be a smooth $H$-fixed point and consider the blow up
$B_xX$ with exceptional divisor $E\cong \p^{n-1}$. 
The $H$-action lifts to $B_xX$ and so we get an $H$-action on $E$
which has a fixed point by (\ref{class.defn}.1).  Since $Y$ is
proper, the induced rational map $B_xX\to X\map Y$ is defined
outside a subset of codimension at least 2. Thus we get an
$H$-equivariant rational map $E\map Y$. By induction, there is a 
fixed point on $Y$.\qed

\begin{rem} If $H$ does not satisfy the  conditions
of (\ref{class.defn})
then (\ref{Going-down}) fails for some actions.
Indeed, let $H\to GL(n,K)$ be a representation without an 
$H$-eigenvector.  This gives an $H$-action on $\p^n$ with a single fixed
point $Q\in \p^n$. The corresponding action on $B_Q\p^n$ has no fixed
points.
\end{rem}

\begin{prop}[Going up]\label{Going-up}
  Let $K$ be an algebraically closed field
and
$H$ a finite abelian group of prime power order $q^n$ 
($q$ is allowed to coincide with $\chr K$). Let $p:X\map Z$ be an
$H$-equivariant map of irreducible
$K$-schemes. Assume that
\begin{enumerate}
\item $f$ is generically finite, dominant and $q\not\vert\deg (X/Z)$,
\item $H$ has a smooth fixed point on $Z$, and
\item $X$ is proper.
\end{enumerate}
\noindent Then $H$ has a fixed point on $X$.
Moreover, if $X\map Y$ is an $H$-equivariant
map to a proper $K$-scheme then $H$ has a fixed point on $Y$.
\end{prop}

Proof.  The proof is by induction on $\dim Z$. 
The  case  $\dim Z=0$ is clear.

Let $z\in Z$ be a smooth fixed point and $E\subset B_zZ$
the exceptional divisor. Let $\bar p:\bar X\to B_zZ$
denote the normalization of $B_zZ$ in the field of rational
functions of $X$ and $F_i\subset \bar X$ the divisors lying over
$E$.  $H$ acts on the set   $\{F_i\}$. Let ${\cal F}_j$
denote the   $H$-orbits and in each pick a divisor
$F^*_j\in {\cal F}_j$.
By the ramification formula,
$$
\deg (X/Z)=\sum_j 
|{\cal F}_j|\cdot\deg (F^*_j/E)\cdot e(\bar p, F^*_j)
$$
 where
$e(\bar p,F^*_j)$ denotes
the ramification index of $\bar p$ at the generic point of $F^*_j$.
 Since $\deg (X/Z)$ is not divisible by $q$,
there is an orbit ${\cal F}_0$ consisting of  a single element
$F^*_0$  such that $\deg (F^*_0/E)$ is not divisible by $q$.

We have $H$-equivariant rational maps 
$F^*_0\map E$, $F^*_0\map Z$ and $F^*_0\map Y$.
By induction $H$ has a fixed point on $Y$, and also on $Z$.
\qed

\begin{rem} We see from the proof that (\ref{Going-up})
also holds  if $H$ is abelian and only one of the prime divisors of
$|H|$ is less than $\deg (X/Z)$.
\end{rem}

The method   also gives a simpler proof of a result
of \cite{nishi}. One can view this as a version of
(\ref{Going-down}) where $H$ is the absolute Galois group of $K$.

\begin{prop}[Nishimura lemma]\label{nishimura}
  Let $K$ be a 
field  and $f:X\map Y$ a 
map of
$K$-schemes. Assume that
\begin{enumerate}
\item $X$ has a smooth $K$-point, and
\item $Y$ is proper.
\end{enumerate}
\noindent Then $Y$ has a   $K$-point.
\end{prop}

Proof. The proof is by induction on $\dim X$. The  case  $\dim
X=0$ is clear.

Let $x\in X$ be a smooth $K$-point and consider the blow up
$B_xX$ with exceptional divisor $E\cong \p^{n-1}$. 
$E$ has smooth $K$-points.
 Since $Y$ is
proper, the induced rational map $B_xX\to X\map Y$ is defined
outside a subset of codimension at least 2 and we get a rational
map $E\map Y$. By induction, there is a  $K$-point on $Y$.\qed

\begin{rem}  One can combine (\ref{Going-down}) and 
(\ref{nishimura}) if we know that any $H$-representation has an
eigenvector defined over $K$. There are two interesting cases
where this condition holds:
\begin{enumerate}
\item  $H$ is Abelian of order $n$ and $K$ contains all $n$th
roots of unity.
\item $H$ is nilpotent and its order is a power of $\chr K$.
\end{enumerate}
\end{rem}

\vskip1cm

\noindent University of Utah, Salt Lake City UT 84112 

\begin{verbatim}kollar@math.utah.edu\end{verbatim}

\noindent Mathematical Institut, Budapest, PO.Box 127, 1364 Hungary

\begin{verbatim}endre@math-inst.hu\end{verbatim}

\end{document}